\documentclass[11pt]{article}
\usepackage{amsmath,amsfonts, amssymb,amsthm,multicol,graphicx,tikz, enumitem,asymptote}
\usepackage{makeidx}
\pdfoutput=1
\renewcommand{\a}{\alpha}

\newcommand{\ds}{\displaystyle}
\newcommand{\ep}{\epsilon}
\newcommand{\PP}{\mathbb{P}}
\newcommand{\bS}{\mathbb{S}}

\newcommand{\cF}{\mathcal{F}}

\newcommand{\cE}{\mathcal{E}}

\newcommand{\cI}{\mathcal{I}}
\newcommand{\cT}{\mathcal{T}}

\newcommand{\EE}{\mathbb{E}}

\newtheorem{prop}{Proposition}   
\newtheorem{lem}[prop]{Lemma} 
\newtheorem{thm}[prop]{Theorem} 
\newtheorem{quest}[prop]{Question}
\newtheorem{conj}[prop]{Conjecture}

\title{Two-Point Concentration of the Independence Number of the Random Graph}
\author{Tom Bohman\thanks{This work was supported in part by a grant from the Simons Foundation (587088, TB)} \  and Jakob Hofstad }

\begin{document}

\maketitle

\begin{center}
    {\bf Abstract}
\end{center}

\noindent
We show that the independence number of $ G_{n,p}$ is concentrated on two values if $ n^{-2/3+ \epsilon} < p \le 1$. This result is roughly best possible as an argument of Sah and Sawhney shows that the independence number is not, in general, concentrated on two values for $ p  = o \left( (\log(n)/n)^{2/3} \right)$. The extent of concentration of the independence number of $ G_{n,p}$ for $ \omega(1/n) <  p \le n^{-2/3}$ remains an interesting open question.

\section{Introduction}

An independent set in a graph $ G = (V,E)$ is a set of vertices that contains no edge of the graph. The independence number
of $G$, denoted $\alpha(G)$, is the maximum number of vertices in an independent set in $G$. The binomial random graph $ G_{n,p}$ has
a vertex set of cardinality $n$, and each potential edge appears independently with probability $p$, where $p=p(n)$ can vary with $n$.
We say that an integer-valued random variable $X$ defined on $G_{n,p}$ is concentrated on
$k = k(n)$ values if there is a function $f(n)$ such that 
\[ \lim_{n \to \infty} \PP( X \in \{f(n), \dots ,f(n)+k-1\} ) =1. \]
In this work we address the following question: for which probabilities $p =p(n)$
is $ \alpha( G_{n,p})$ concentrated on two values? 

The independence number of $ G_{n,p}$ has been a central issue in the study of the random graph since the beginning. In the 1970's, 
Bollob\'as and Erd\H{o}s \cite{1976cliques} and Matula \cite{matula} independently showed that $ \alpha( G_{n,p})$ is concentrated on 
two values when $p$
is a constant. In the late 1980's, Frieze \cite{frieze} showed that if $ \omega(1/n) < p < o(1)$ then
\begin{equation} \label{eq:history}
\a( G_{n,p}) = \frac{2}{p} \left[ \log(np) - \log\log(np) + \log(e/2) \pm o(1) \right]  
\end{equation}
with high probability. This celebrated result combines the second moment method with a large
deviation inequality in a surprising way. Dani and Moore \cite{dm} give the best known lower bound on the independence
number of $ G_{n,p}$ when $p$ is a constant over $n$. For further discussion of $ \a(G_{n,p})$ see the canonical random graphs texts \cite{bela} 
\cite{alanmichal} \cite{JLR}.

The algorithmic question of finding a large independent set in a random graph is also a central issue. An influential
work of Karp and Sipser \cite{ks} introduced a simple randomized algorithm that determines the independence number of $ G_{n,p}$ for
$ p = c/n$ such that $c < e$.  However, the computational problem of finding a maximum independent set in $ G_{n,p}$ for larger $p$ appears to be a very difficult problem; for example, Coja-Oghlan and Efthymiou \cite{amin} showed that the space of independent sets in $ G_{n,p}$ that are larger than roughly half the independence number form an `intricately ragged landscape' that foils local search algorithms.

Interest in the concentration of the independence number of $ G_{n,p}$ is also motivated by the study of the 
chromatic number of $ G_{n,p}$. An influential series of works starting with Shamir and Spencer \cite{ss} and culminating
in the result of Alon and Krivelevich \cite{ak} establish two point concentration of the chromatic number of $ G_{n,p}$
for $ p < n^{-1/2 - \epsilon}$  (also see \cite{luczak}, \cite{an}, \cite{cps}, \cite{lutz}). More recently Heckel \cite{h} showed that
the chromatic number of $ G_{n,1/2}$ is
{\bf not} concentrated on any interval of length smaller then $ n^{1/4- \epsilon}$. Heckel and Riordan \cite{hr} then 
improved this lower bound on the length of the interval on which the chromatic number is concentrated to $ n^{1/2 - \epsilon}$. They also state a fascinating conjecture regarding fine detail of the concentration of the chromatic number of $ G_{n,1/2}$.
Determining the extent of concentration
of the chromatic number of $ G_{n,p}$ for $ n^{-1/2} \le p < o(1)$ is another central and interesting question.
Two point concentration of the domination number of $ G_{n,p}$ has also been established for a wide range of values of $p$ \cite{gls}.

Despite this extensive interest in two point concentration of closely related parameters, two point concentration of
$ \alpha(G_{n,p})$ itself has not been addressed since the early pioneering results of Bollob\'as and Erd\H{o}s and Matula. Our main result is as follows.
\begin{thm} \label{MainLoose} Let $ \epsilon >0$.
If $ n^{-2/3+\ep} < p \le 1$ then $ \alpha( G_{n,p})$ is concentrated on two values. 
\end{thm}

Soon after the original version of this manuscript was posted to the ArXiv, Sah and Sawhney observed that 
Theorem~\ref{MainLoose} is roughly best possible. 
\begin{thm}[Sah and Sawhney \cite{sahs}]\label{SS}
Suppose $ p = p(n)$ satisfies $ p = \omega(1/n)$ and  $ p <  (\log(n)/(12n))^{2/3}$. Then there exists $ q = q(n)$ such that $ p \le q \le 2p$ and $ \alpha( G_{n,q})$ is not concentrated on 
fewer than $ n^{-1}p^{-3/2} \log(np) / 2$ values. 
\end{thm}
\noindent
Note that this anti-concentration statement is somewhat weak as it is compatible with two point concentration for some functions $ p < n^{-2/3}$. On the other hand, it is natural to suspect that the extent of concentration of the independence number is monotone in the regime. See the Conclusion for further discussion.

Now, the second moment method is the main tool in the proof of Theorem~\ref{MainLoose}. 
Let $X_k$ be the number of independent
sets of size $k$ in $G_{n,p}$. We have
\[ \EE[ X_k ] = \binom{n}{k} (1-p)^{\binom{k}{2}}, \]
and if $k$ is large enough for this expectation to vanish then there are no independent sets of size $k$ with high probability. This simple observation provides
a general upper bound on $ \alpha( G_{n,p})$. As far as we are aware, this is the only upper bound on $ \alpha( G_{n,p})$ that appears in the
literature.  It turns out this upper bound is sufficient to establish two point concentration when 
$ n^{-1/2 + \epsilon} < p \le 1$.
(This was recently observed by the authors \cite{bh} for $n^{-1/3+ \epsilon} < p < o(1)$.)
However, for smaller values of $p$ this upper bound turns out not to be sufficient for two point concentration 
as the independence number is actually somewhat
smaller than the largest value of $k$ for which $\EE[X_k]= \omega(1)$ with high probability. 

In the regime
$ n^{-2/3+ \epsilon} < p < n^{-1/2 + \epsilon}$ we need to consider a more general object.  We define an {\em augmented independent set of order $k$}
to be a set $S$ of $k+r$ vertices such that the graph induced on $S$ is a matching with $r$ edges and every vertex $v \notin S$ has at least
two neighbors in $S$. Clearly, such an augmented independent set contains $2^r$ independent sets of size $k$, and so $\EE[X_k]$ can be quite large even when the
expected number of augmented independent sets of order $k$ vanishes. As the appearance of a single independent set of size $k$ implies the appearance of
an augmented independent set of order at least $k$, it follows from this observation that the appearance of augmented independent sets containing large induced matchings can create a situation in which $ \alpha( G_{n,p}) $ is less then $k$ whp even
though $ \EE[X_k]$ tends to infinity. Furthermore, it turns out the variance in the number of augmented independent sets of order $k$ is small enough
to allow us to establish two point concentration using the second moment method in the regime in question.

The remainder of this paper is organized as follows. In the following section we make some preliminary calculations. 
We then prove our main result, Theorem~\ref{MainLoose}, in Section~3. 
Theorem~\ref{SS} is proved in Section~\ref{sec:ss}. 
In Section~\ref{sec:anti} we observe that 
if $ p = c/n$ for some constant $c$ then there is a constant $K$ such that 
$ \alpha( G_{n,p})$ is not concentrated on any interval of 
length $ K\sqrt{n}$.
We conclude by making some remarks regarding the concentration $ \a ( G_{n,p})$ for
$ \omega(1/n) < p \le n^{-2/3}$ in Section~\ref{sec:concl}.

We mentioned above that
the second moment method applied to $X_k$ is sufficient to establish two point concentration down to $ p = n^{-1/3+ \epsilon}$ (see \cite{bh} for this proof). 
It turns out that
the second moment method applied to the number of maximal independent sets of size $k$ suffices to establish two point concentration down to 
$ p = n^{-1/2 + \epsilon}$. We 
include this proof in an Appendix.

\section{Preliminary Calculations}

Throughout Sections~2~and~3 we restrict our attention to $p = p(n)$ where $p =n^{-2/3 + \ep} < p < 1/ \log(n)^2$ for some $\ep>0$. This is valid as the
classical results of Bollob\'as and Erd\H{o}s and Matula treat the case $ p > 1/ \log(n)^2$.  It is a classical result
of Frieze that for all such $p$ we have
\begin{align}
    \alpha( G_{n,p}) = \frac{2 \left[ \log(np) - \log(\log(np)) + \log(e/2) \pm o(1) \right] }{p}
\end{align}
with high probability. As we are interested in refinements of this result we will need some estimates regarding
numbers in this range.  

Recall that $ X_k$ is the number of independent sets of size $k$ in $ G_{n,p}$.
The upper bound in Frieze's result is achieved by simply identifying a value of $k$ such that $\EE [X_k] = o(1)$. 
In order to make comparisons with this bound, we define $ k_x$ to be the largest value of $k$ for which
\[ \EE[X_k] = \binom{n}{k} (1-p)^{\binom{k}{2}} > n^{2\epsilon}. \] 
We note in passing that the exact value in the exponent of this cut-off function is not crucial; we simply need it to be small (in particular, we need $ \EE[ X_{k_x+2}] = o(1)$) and greater than the exponent in the polynomial function bound that we use to define the variable $k_z$ below.

Now suppose $ k = k_x \pm o(1/p)$.  
In other words, consider 
\begin{equation}
k = \frac{2 \left[ \log(np) - \log(\log(np))   + \log( e/2) \pm o(1) \right] }{p}.   \label{approx}  
\end{equation}
(Note that impact of the $ 2 \epsilon$ in the definition of $ k_x$ is absorbed in the $ o(1)$ term.)
Our first observation is that, for such a $k$, we have
\begin{multline*}
\frac{ne}k (1-p)^{k/2} = (1+o(1)) \frac{ne}k \exp\{ - \log(np) + \log( \log(np)) - \log( e/2) + o(1)\} \\ = ( 1 +o(1)) \frac{ \log(np) }{ \log(np) - \log\log(np)  + \log(e/2) - o(1)}   
= 1+o(1).
\end{multline*}
So, we have
\begin{align} 
 (1-p)^{k/2} = (1+o(1)) \frac{k}{e n}. \label{k pre-approx}
\end{align}
It follows that we have
\begin{gather} \label{first max factor}
\left( 1 - ( 1 -p)^k\right)^{n-k} =  \left( 1 - (1+o(1)) \left( \frac{k}{en}\right)^2 \right)^{n-k}
= e^{ - (1 +o(1)) \frac{k^2}{ e^2 n}}. 
\end{gather}
Similarly,
\begin{equation}
\label{eq:maxxy}
 \left[ 1 - (1-p)^k - pk(1-p)^{k-1} \right] ^{n-k} = e^{ - (1 +o(1)) \frac{ p k^3}{ e^2 n} }.
 \end{equation}

\section{Two-point concentration}

In this section we show that two point concentration of $\a(G_{n,p})$ persists for $p$ as low as $n^{-2/3+\ep}$, but at a value that does not equal $k_x$ for small $p$. Our main technical result is as follows. We write $ f(n) \sim g(n)$ if $ \lim_{n \to \infty} f/g =1 $.
\begin{thm} \label{Main}
If $ n^{-2/3+\ep} < p < 1 / \log(n)^2 $ for some $\ep>0$, then there exists an integer $k_{z} = k_{z}(n)$ such that $\a(G_{n,p}) \in \{k_{z},k_{z}+1\}$ whp. Furthermore, we have 
\[ k_x - k_{z} \begin{cases} = 0 & \text{ if } p = \omega\left( \log(n) n^{-1/2} \right)\\
=  \xi_n & \text{ if } p = C \log(n) n^{-1/2} \\
\sim  \frac{k_x^2}{e^2 n} \sim \frac{ 4 \log(np)^2}{ e^2 p^2 n} & \text{ if } p = o\left( \log(n) n^{-1/2} \right) \end{cases} \]
where $ \xi_n \in \left( \frac{1}{ e^2C^2}  -\frac{5}{2}, \frac{1}{ e^2C^2} +\frac{3}{2} \right) $. 

\end{thm}
\noindent
We emphasize that the sequence $ \xi_n$ for the case 
$ p = C \log(n) n^{-1/2}$ does not converge to a particular value; rather, it ranges over 
the specified interval. Indeed, $ \xi_n$ is more precisely given by the following expression
\[ \xi_n =  \left\lceil  \frac{1}{e^2C^2} - \frac{2\log(\EE[X_{k_x}])}{\log(n)} + 2\epsilon \pm o(1)   \right\rceil. \]
As $ \EE[ X_{k_x}]$ can take a wide range of values as $n$ varies, we have persistent variation of $ \xi_n$ over a small list of values.
The remainder of this Section is a proof of Theorem~\ref{Main}.

Recall that an
{\it augmented independent set of order $k$} in a graph $G$ is defined to be a set of vertices $S$ for which there exists some $r \geq 0$ such that $|S| = k+r$ and $S$ contains exactly $r$ edges that form a matching (i.e. these edges are pairwise disjoint) and all vertices outside of $S$ have at least 2 neighbors in $S$. To motivate this definition, first note that such a set $S$ contains $2^r$ independent sets of size $k$, so this structure is well suited for isolating large variations in the number of independent sets of size $k$. The fact that we are interested in studying the number of independent sets of size $k$ where $k$ is close to the independence number of $G$ motivates the condition regarding vertices outside of $S$. Indeed, suppose $S$ is a set such that the induced graph on $S$ is a (not necessarily perfect) matching and $v$ is a vertex outside of $S$ that has one neighbor in $S$. If $v$ is adjacent to a vertex $u$ that is not in the matching then adding $v$ to $S$ creates an augmented independent set of order $k$ with an additional edge in its matching. So in this situation we would want to include $v$ in $S$ in order to isolate as much variation in the number of independent sets of $k$ as possible. On the other hand, if $v$ is adjacent to a vertex $u$ that is part of the matching then $ S \setminus \{u\} \cup \{v\}$ 
contains an independent set of size $k+1$.

\begin{lem} \label{MNI connection}
For any graph $G$, $\a(G) = k$ if and only if $G$ has an augmented independent set of order $k$, but no augmented independent set of any larger order.
\end{lem}
\begin{proof}
Let $ \hat{\alpha}(G)$ be the largest $k$ for which there is an augmented independent set of order $k$. As an augmented independent set of order $k$ contains an
independent set on $k$ vertices we have
\[  \hat{\alpha}(G) \le \alpha(G).   \]
Now suppose that $G$ has a maximum independent set $S$ of size $k$. Let $T$ be a maximum set of vertices that contains $S$ and has the additional property
that the graph induced on $T$ is a matching. We claim that $T$ is an augmented independent set; that is, every vertex outside of $T$ has 2 neighbors in $T$.
Clearly every vertex outside $T$ has a neighbor in $T$ (since $S$ is a maximum independent set).  If $v$ is a vertex outside
of $T$ with only one neighbor in $T$ then this neighbor $u$ is in one of the matching edges (since $T$ is chosen to be a maximum set with the given property), but in this situation
the set $ (T \setminus \{u\}) \cup \{v\}$ contains an independent set that is larger than $S$, which is a contradiction. So, $T$ is an augmented independent set of order at least $k$. Therefore,
\[ \alpha(G) \le \hat{\alpha}(G). \]

\end{proof}

We prove Theorem~\ref{Main} by applying the second moment method to the random variable $Z$, which counts the number of augmented independent sets of order $k_z$ with a prescribed number of edges $r_z$. This requires some relatively delicate estimates to determine the optimal values of the parameters $k_z$ and $r_z$. We begin with those calculations. We then move on to the first moment and second moment calculations for $Z$ that suffice to prove Theorem~\ref{Main}.

\subsection{Defining and estimating $k_z$ and $r_z$}

In this section we define $ k_z$ and $r_z$, and establish a number of estimates of these quantities.

To start,
let $E(n,k,r)$ denote the expected number of augmented independent sets of order $k$ with exactly $r$ edges in $G_{n,p}$. There are $\binom{n}{k+r}$ ways to choose the location of our set, and given this set there are $\frac{(k+r)!}{(k-r)! 2^r r!}$ ways to choose the locations of the $r$ disjoint edges. The probability of choosing internal edges accordingly is \[p^r(1-p)^{\binom{k+r}{2}-r},\] and finally, the probability that every vertex outside our set has at least two neighbors within the set is $ F(n,k,r)^{n-k-r}$ where we set
\begin{align}
    F(n,k,r) := 1 - (1-p)^{k+r} - (k+r)p(1-p)^{k+r-1}. \label{F}
\end{align}
Thus, we have
\begin{align}
    E(n,k,r) 
    =\binom{n}{k+r}\frac{(k+r)!}{(k-r)! 2^r r!}p^r(1-p)^{\binom{k+r}{2}-r}  F(n,k,r)^{n-k-r} . \label{first moment formula}
\end{align}
Now we define $k_z = k_z(n)$ to be the largest $k \le k_x$ for which there exists an $r$ such that
$ E(n,k,r) >  n^{\epsilon}$.\label{footnote 1}\footnote{We impose the condition $ k \le k_x$ to avoid situations in which $k_z > k_x$. In particular, there are certain values of $n$ and $p$ for which we would have $k_z = k_x + 1$. For example, this would happen if $n^{3\epsilon/2} < \EE[X_{k_x+1}] <n^{2 \epsilon}$.}

Our first observation is that while $k_z$ is not equal to $k_x$ when $p$ is sufficiently small, the difference between the
two numbers is relatively small.
\begin{lem} \label{lem:zbig}
If $ p = \omega( \log(n)/ \sqrt{n} )$ then $ k_z = k_x$. Otherwise we have
\begin{equation*}
    k_z \ge k_x - O \left( \frac{ k_x^2}{n} \right) \ge k_x - o( k_x^{1/2}). %
\end{equation*}
\end{lem}
\begin{proof}
This observation follows from considering $r=0$ and $ k = k_x - \frac{ C k_x^2}{n}$ where $C$ is a sufficiently large constant and applying (\ref{k pre-approx}) and (\ref{eq:maxxy}).  
We first note that if $ p = \omega( \log(n)/ \sqrt{n})$ then the expression in (\ref{eq:maxxy}) is subpolynomial. In this case we have $ E(n,k_x,0) > n^{\epsilon}$ and $ k_z = k_x$.

Restricting our attention to $ p = O ( \log(n)/ \sqrt{n})$, note that
for $ k = k_x - O( k_x^2/n)$ we have 
\begin{equation} \label{change in binomial}
    \frac{ \binom{n}{k} (1-p)^{\binom{k}{2}} }{ \binom{n}{k+1} (1-p)^{\binom{k+1}{2}}  } = \frac{ k+1}{n-k} (1-p)^{-k} = (1+o(1)) \frac{ ne^2}{k},
\end{equation}
and thus
\begin{equation*}
    \begin{split}
    E(n,k,0) &=\binom{n}{k} (1-p)^{ \binom{k}{2} } \exp \left\{- (1+o(1)) \frac{p k^3}{ e^2 n} \right\} \nonumber \\& \ge
 \left( \frac{ ne^2}{ 2k} \right)^{ \frac{ C k^2}{n}} \exp \left\{- (1+o(1)) \frac{p k^3}{ e^2 n} \right\}.
\end{split}
\end{equation*}
Therefore $E(n,k,0) > n^{\epsilon}$ if $C$ is a sufficiently large constant.
\end{proof}
\noindent
Note that we may henceforth assume that the estimates given in (\ref{approx})-(\ref{eq:maxxy}) hold for $k_z \le k \le k_x$.

Next, for $ k_z \le k \le k_x$  define \[r_M(k) = r_M(n,k) := \arg \ds\max_{r} \{E(n,k,r)\}. \]
We are now ready to define the random variable $Z$. \label{define Z} Set $ r_z = r_M(k_z)$ and let the variable $Z$ count the number of augmented independent sets of order $k_z$ with exactly $r_z$ internal edges.  Note that such sets have $k_{z}+r_z$ vertices and  $\EE(Z) = E(n,k_z,r_z) > n^\epsilon$. 

We now move on to more accurate estimates for $k_z$ and $r_z$, starting with an estimate for $r_z$ that we obtain by estimating the ratio of consecutive terms of (\ref{first moment formula}). Note that we always have $r \leq k$, which we will use a couple of times below. First, we show that $F(n,k,r)^{n-k-r}$ will have an insignificant effect on the ratio. One can easily verify, using (\ref{k pre-approx}), that $F(n,k,r) = 1-o(1)$. Furthermore,
\[ F(n,k,r+1) - F(n,k,r) =(k+r)p^2(1 -p)^{k+r-1}\] exactly.
Therefore, we have
\begin{align*}
    \frac{F(n,k,r+1)^{n-k-r-1}}{F(n,k,r)^{n-k-r}} & = \left(\frac{F(n,k,r+1)}{F(n,k,r)}\right)^{n-k-r} F(n,k,r+1)^{-1} 
    \\ & = (1 +o(1)) \left( 1 + \frac{ (k+r)p^2(1 -p)^{k+r-1}}{ F(n,k,r)}  \right)^{n-k-r} \\
    & = (1 +o(1)) \left( 1 + (1+o(1)) (k+r)p^2(1 -p)^{k+r-1}  \right)^{n-k-r}.
\end{align*}
Since $$(k+r)p^2(1-p)^{k+r-1}(n-k-r) \leq 2 k p^2 n(1-p)^{k-1} < \frac{k^3 p^2}{n} = o(1), $$ it follows that 
\begin{equation}
    \label{eq:ratiosmall}
\frac{F(n,k,r+1)^{n-k-r-1}}{F(n,k,r)^{n-k-r}} = 1+o(1), 
\end{equation}
as desired.

With this bound on the ratio $ F(n,k,r+1)^{n-k-r-1}/ F(n,k,r)^{n-k-r}$ in hand, and applying (\ref{k pre-approx})
as in (\ref{change in binomial}),
we conclude that
if $r = o(1/p)$ then we have
\begin{align}
    \frac{E(n,k,r+1)}{E(n,k,r)} \nonumber &=(1+o(1)) 
    \frac{ \binom{n}{k+r+1} (1-p)^{\binom{k+r+1}{2}} }{ \binom{n}{k+r} (1-p)^{\binom{k+r}{2}}  } \left(\frac{(k+r+1)(k-r)p}{2(r+1)}\right) \\ 
    & = (1+o(1)) \left(\frac{k}{e^2 n}\right) \left(\frac{(k+r+1)(k-r)p}{2(r+1)}\right)\\&=
    (1+o(1))\frac{k^3p}{2e^2n(r+1)}(1-r^2/k^2). \label{ratio simplified}
\end{align}
This calculation provides the desired estimate for $r_z$.
\begin{lem} \label{rz without k}
\[ r_z = \left\lceil (1+o(1)) \frac{4 \log(np)^3 }{ e^2 n p^2} -1  \right\rceil. \] 
\end{lem}
\begin{proof}
First note that the form of (\ref{ratio simplified}) suggests that $ r_M(k)$ should
be roughly $ k^3p/n$ as this is the regime in which the ratio is roughly one. 
Since $k^3p/n = o(1/p)$, the estimate given by (\ref{ratio simplified}) holds
in the vicinity of $ r = k^3p/n$. Furthermore, for $  r= \Omega(1/p)$ the arguments
above can be easily adapted to show that
$ E(n,k,r+1)/E(n,k,r)$ is bounded above by the expression in (\ref{ratio simplified}).
It follows that we have
\begin{equation}
r \ge k^3p/(e^2n) \ \ \ \Rightarrow \ \ \ \frac{E(n,k,r+1)}{E(n,k,r)} < 1.
    \label{ratio simplified modified}
\end{equation}
Hence, we can conclude that
\begin{equation}
\label{eq:looseupper}
r_M(k) \le \frac{k^3p}{e^2n} = O\left( \frac{\log(np) k^2}{n} \right). 
\end{equation}
Note that if $r = O( k^3p/n )$ then the $(1-r^2/k^2) $ term in (\ref{ratio simplified}) is equal to $ 1+o(1)$.
We conclude that we have 
\begin{equation}
\label{eq:rmprecise}
r_M(k) = \left\lceil \frac{k^3p}{ 2 e^2n} (1+o(1)) -1  \right\rceil. 
\end{equation}
In light of Lemma~\ref{lem:zbig}, the proof is complete.
\end{proof}
We note in passing that Lemma~\ref{rz without k} implies that $ r_z $ becomes relevant at 
$ p = \Theta( \log(n)^{3/2}/ n^{1/2}) $, as $ r_z=0$ for $p$ larger than this regime.  
Furthermore, for all $p$ that we contemplate in this section (i.e. for $p> n^{-2/3 + \ep}$) 
we have 
\begin{equation}
    \label{eq:rbound}
    r_z = o( k_x^{1/2} ).
\end{equation}

Finally we turn to estimating $ k_z $. Note that we may assume $ p = O( \log(n) n^{-1/2})$ as we have
$ k_z= k_x $ for larger $p$ by Lemma~\ref{lem:zbig}. Note further that (\ref{eq:rmprecise}) implies that $ r_M(k) = \Omega( \log(n))$ for $ k_z \le k \le k_x$ for $p$ in this
regime. Recall that, appealing to (\ref{eq:maxxy}), we have \[ E(n,k_x,0) = \EE[ X_{k_x}]  e^{- (1+o(1)) \frac{p k_x^3}{ e^2 n} }. \]
Next observe that, since $F(n,k+1,r) = F(n,k,r+1)$, we can easily adapt 
the argument that proves (\ref{eq:ratiosmall}) to establish
\[ \frac{F(n,k+1,r)^{n-k-r-1}}{F(n,k,r)^{n-k-r}} = 1+o(1) \]
for $ k_z \le k \le k_x$ and $r = o(1/p)$. 
So we
can calculate the ratio over a change in $k$ in the same way that we established (\ref{ratio simplified}). Thus, for $ k_z \le k \le k_x$ and $r = o(1/p)$ we have 
\begin{align}
    \frac{E(n,k+1,r)}{E(n,k,r)} = \frac{k}{e^2 n}(1+o(1)) \label{ratio over k}.
\end{align}
Using (\ref{eq:maxxy}), (\ref{ratio over k}), (\ref{ratio simplified}), and (\ref{eq:rmprecise}), we write
\begin{equation*}
\begin{split}
    E(n,k,r_M(k)) & = E(n,k_x,0) \prod_{\ell=k}^{k_x-1} \frac{E (n,\ell,0)}{ E(n, \ell+1,0) } \prod_{r=0}^{r_M(k)-1} \frac{E(n,k,r+1)}{ E(n,k,r)} \\
    & = \EE[ X_{k_x}]  e^{- (1+o(1)) \frac{p k_x^3}{ e^2 n} } \left( \frac{e^2n}{k} (1+o(1))  \right)^{k_x-k}   \left((1+o(1))\frac{k^3p}{2e^2n} \right)^{r_M(k)} \frac{1}{ r_M(k)!} \\
    & =  \EE[ X_{k_x}]  e^{- (1+o(1)) \frac{p k_x^3}{ 2 e^2 n} } \left( \frac{e^2n}{k} (1+o(1))  \right)^{k_x-k}. 
    \end{split}
\end{equation*}
This estimate suffices to establish the following.
\begin{lem}\label{lem:kz} 
\begin{align*}
n^{-2/3 + \epsilon} < p \ll \log(n) n^{-1/2} &  \hskip5mm \Rightarrow &
k_x - k_z = (1+o(1)) \frac{ k_x^2 }{ e^2 n  } = (1+o(1)) \frac{ 4 \log(np)^2 }{ e^2p^2n} \\
p = C \log(n) n^{-1/2} &  \hskip5mm \Rightarrow &   k_x - k_z =  \left\lceil \frac{1}{e^2C^2} - \frac{2\log(\EE[X_{k_x}])}{\log(n)} + 2 \epsilon \pm o(1)   \right\rceil.
\end{align*}
\end{lem}
\begin{proof}
First note that for $p$ in these ranges we have
\[ \log( n/k) = (1+o(1)) \log(np) = (1+o(1))\frac{pk}{2}.\]
It follows that $ k_x - k_z$ is the smallest integer $ \kappa$ such that
\[ \log \EE[ X_{k}] - ( 1+o(1)) \frac{ k_x^3p}{2e^2n} + \kappa (1+o(1))\frac{pk_x}{2} > \epsilon \log(n).  \]
If $ p \ll \log(n) n^{-1/2}$ then $ \log \EE[X_k] = O( \log(n)) = o( k_x^3p/n)$ and the first part of the Lemma follows. On the other hand, if $ p = C \log(n) n^{-1/2}$ then  $ k_x^3 p /n = (1+o(1))\log(n)/C^2$ and 
$ p k_x = \log(n) (1+o(1))$ and we recover the second part of the Lemma.

\end{proof}

\subsection{First moment} \label{sec:firstmoment}

Here we show that, with high probability, no augmented independent set of order $k$ appears for any $ k_{z}+2 \le k \le k_x+1$. This is sufficient because we can simply consider the
expected number of independent sets of size $k_x + 2$ to rule out the appearance of any larger independent set. 
We emphasize that we are restricting our attention to $k$ that satisfy (\ref{approx}). We begin with a simple observation.
\begin{lem} \label{bulk}
For all $k$ that satisfy (\ref{approx}) we have 
$$ \ds\sum_{r=0}^{k} E(n,k,r) \leq 4 \sum_{r=0}^{2r_M(k)} E(n,k,r).$$
\end{lem}
\begin{proof}
We begin by noting that (\ref{ratio simplified}) can be adapted to show
\[ \frac{ E(n,k,r+1) }{ E(n,k,r)}  \le (1+o(1))\frac{k^3p}{2e^2n(r+1)}. \]
Now we consider cases depending on the value of $ r_M = r_M(k)$.  
First note that if $ r_M = 0$ then we have $ k^3p/ (2e^2n) < 1 +o(1) $ and
\[ \sum_{r=0}^k E(n,k,r) \le (1+o(1)) E(n,k,0) \sum_{r=0}^\infty \frac{1}{r!} \le 4 E(n,k,0). \] 
Now suppose $ r_M \ge 1$.  Recalling (\ref{eq:rmprecise}), note that if $r > 2r_M$ then we have \begin{equation} E(n,k,r+1)/E(n,k,r) \le 2/3.  \label{eq:babystep}
\end{equation}
So we have
\begin{multline*}
\ds\sum_{r=0}^{k} E(n,k,r) =  \sum_{r=0}^{2r_M} E(n,k,r) + \sum_{r =2r_M+1}^k E(n,k,r) \\ \le  
\left(\sum_{r=0}^{2r_M} E(n,k,r)\right) + 3 E(n,k,2r_M+1)
\le \left(\sum_{r=0}^{2r_M} E(n,k,r)\right) + 3 E(n,k,2r_M). 
\end{multline*}
\end{proof}

Note that it follows from the proof of Lemma~\ref{bulk} that we have
\[ E(n, k+1,k+1) \leq \frac{2}{3} E(n,k+1,k). \]
Applying this observation together with Lemma~\ref{bulk}, and (\ref{ratio over k}), we have
\begin{align*}
    \sum_{r=0}^{k+1}E(n,k+1,r) &< \frac{5}{3} \sum_{r=0}^k E(n,k+1,r) \\&<
    \frac{20}{3} \sum_{r=0}^{2 r_M(k+1)} E(n,k+1,r)\\&<
    \frac{k}{n} \sum_{r=0}^{2 r_M(k+1)} E(n,k,r)\\&<
    \frac{k}{n} \sum_{r=0}^{k}E(n,k,r).
\end{align*}
Therefore, by Lemma~\ref{bulk} and (\ref{eq:looseupper}), we have
\begin{align*}
    \ds\sum_{k=k_z+2}^{k_x+1} \left(\sum_{r=0}^{k}E(n,k,r)\right) &< \frac{2k_z}{n} \sum_{r=0}^{k_z+1} E(n,k_z+1,r) \\&\leq
    \frac{8k_z}{n} \sum_{r=0}^{2r_M(k_z+1)} E(n,k_z+1,r) \\&\leq
    \left(\frac{8k_z}{n}\right)\left(\frac{3k_z^3 p}{e^2 n}\right) E(n,k_z+1,r_M( k_z+1)) \\&\leq
    \left(\frac{8k_z}{n}\right)\left(\frac{3k_z^3 p}{e^2 n}\right) n^\ep \\&=
    O(\log(n)^4n^{-2\ep}).
\end{align*}

Hence, whp no augmented independent set of order $k$ for any $k \geq k_z+2$ appears, and by Lemma~\ref{MNI connection}, we have $\a(G) \leq k_z+1$ whp.

\subsection{Second moment method applied to $Z$}

Here we prove by the second moment method that an augmented independent set of order $k_z$ appears with high probability; more specifically, we prove that such a set with $k_z+r_z$ vertices and $r_z = r_M( k_z)$ internal edges appears whp. For convenience we let 
\[k = k_z, \ \ \ \ \ \ r = r_z = r_M(k_z), \ \ \ \text{and} \ \ \ \Tilde{k} = k + r\] throughout the rest of this subsection. Recall that the random variable $Z$ counts the number of augmented independent sets of order $k$ with $r$ edges, and by the definition of $r_M$ and $Z$ we have
\[ \EE(Z) > n^\ep,\]
and our goal is to show that $ Z>0$ whp.

We now break up $Z$ into a sum of indicator random variables.
Let $ \bS$ be the collection of all pairs $ (S, m_S)$ where $S$ is a set of $\Tilde{k}$ 
vertices and $ m_S$ is a matching consisting of $r$ edges all of which are 
contained in $S$.  Note that
\[ | \bS| =  \binom{n}{\Tilde{k}} \frac{ \Tilde{k}!}{ \left(\Tilde{k}-2r\right)! 2^r r!}.  \]
For each pair $ (S, m_S) \in \bS$, let $ Z_{S,m_S}$ be the indicator random variable for the event that $S$ is an augmented independent set with matching $m_S$.  
We have
\begin{multline*}
    \text{Var}(Z) \leq \sum_{(S,m_S) \in \bS}  \text{Var}(Z_{S,m_S})+ \sum_{S \neq T} \text{Cov}(Z_{S,m_S},Z_{T,m_T}) \\ \leq \EE(Z) + \sum_{S \neq T} \EE(Z_{S,m_S} Z_{T,m_T}) - \EE(Z_{S,m_S})\EE(Z_{T,m_T}),
\end{multline*}
where here and throughout this section such summations are over all $ (S,m_S), (T,m_T) \in \bS$ that have the specified $S$ and/or $T$.
Note that we are making use of the full covariance; we will see below that this is necessary to handle sets $S,T$ 
such that $ |S\cap T|$ is roughly $ k^2/n$. Next, let
\begin{align}
    h_i &:= \frac{1}{\EE(Z)^2}\sum_{|S \cap T| = i} \EE(Z_{S,m_S} Z_{T,m_T}) - \EE(Z_{S,m_S}) \EE(Z_{T,m_T})  \label{h''}    \\&=
    \left(\frac{1 }{ |\bS| }\right)^2\sum_{|S \cap T| = i} \left(\frac{\EE(Z_{T,m_T} | Z_{S, m_S} = 1)}{\EE(Z_{S,m_S})} - 1\right) \label{h'''}.
\end{align}
Applying the second moment method we have $ \PP( Z = 0) \le \text{Var}[Z]/ \EE[Z]^2 $. Thus our goal is to show
\begin{equation*}
\label{eq:TS}
    \sum_{i=0}^{\Tilde{k}-1} h_i = o(1).
\end{equation*}
We emphasize that $ \EE(Z)$ can be significantly different from $ \EE(X)$ here. We consider three ranges of $i$ to establish the desired bounds on $ \sum h_i$.

\noindent
\subsubsection{{\bf Case 1:} $i < \frac{1}{ \log(n) p}$}

Fix a pair $(S, m_S) \in \bS$ and define $ \cE $ to be the event that $ Z_{S,m_S} =1 $.  In other words, $ \cE $ is the event that $S$ gives an augmented independent set with the prescribed matching $m_S$. 
Let $ \cE_T$ be the event $ Z_{T,m_T} =1 $.  Recalling (\ref{h'''}), we seek to bound
\begin{equation}
\label{eq:target}
h_i = \frac{1}{|\bS|} \sum_{T: |S \cap T|=i } \left(  \frac{ \PP( \cE_T \mid \cE)}{ \PP ( \cE) } - 1 \right). 
\end{equation}
Now, for $ \ell = 0, \dots, r $ let $ \bS_\ell$ be the collection of pairs $(T, m_T) \in \bS$ with the property that $ m_S$ and $m_T$ have exactly $\ell$ edges in common.
Define
\[B = \left\lceil \frac{1}{\sqrt{p\ln(n)}}\right\rceil,\]
and, for $ B < i < B^2$, let $\bS_{i,\ell}$ be 
the collection of pairs $(T, m_T) \in \bS_\ell$ such that $|T \cap S| = i$.
Finally, let $\bS'_\ell$ be the collection of pairs $(T, m_T) \in \bS_\ell$ with the
additional property that $ |T \cap S| \le B$. 
We make two claims.

\begin{lem} \label{lem:conditioning}
If $(T, m_T) \in \bS_{i,\ell}$ and $ i < \frac{1}{\log(n)p }$ then 
\[ \PP( \cE_T \mid \cE ) \le e^{O(k^3p^2i/n) +o(1)} \frac{ \PP( \cE) }{ p^\ell (1-p)^{\binom{i}{2}}}. \]
\end{lem}

\begin{proof}

In order to discuss the conditioning on the event $ \cE$, we need to define an additional parameter.  Let $w$ be the number
of edges in $m_T$ that have exactly one vertex in $ S\cap T$ and
let $W$ be the vertices in $T \setminus S$ that are included in such edges. To
be precise we define
\[ W =\{ v \in T \setminus S: \exists \ u \in S \cap T \text{ such that } uv \in m_T   \}. \]
For each vertex $ v \in W$ let $ vv'$ be the edge of $m_T$ that contains $v$.

We condition as follows. First, we simply condition on
all pairs within $S$ appearing as edges and non-edges as required for the event $ \cE$. But we proceed more carefully for vertices $v$ that are not in $S$. For
such vertices we must reveal information about the neighbors of $v$ in $S$. 
We do this by observing whether or not
the potential edges of the form $vu$ where $u \in S$ actually are 
edges one at a time, starting with potential neighbors $ u \in S \setminus T$, and 
stopping as soon
as we observe at least $ 2$ such neighbors. For vertices $v\in W$ we add the restriction
that the edge $ vv' \in m_T $ is the last edge observed in this process.
So, for a vertex $v \in T \setminus S$, this process either observes two edges between $v$ and $ S \setminus T$ - and therefore observes no potential edges between $v$ and $S \cap T$ - or the process observes at most one edge between $v$ and $ S \setminus T$ and therefore observes edges between $v$ and $S \cap T$ and reveals that at least one such edge appears.
In this latter case the conditional probability of $ \cE_T$ is zero {\em unless}
this vertex $v$ is in $W$ and the observed edge is $vv'$.

For each set $ U \subset T \setminus S$ we define $ \cF_{U}$ to be the event that $\cE$ holds and
$U$ is the set of vertices $ v \in T \setminus S$ with the property that our process reveals edges 
between $ v$ and $ S \cap T$. Note that we have
\[ \cE = \bigvee_{ U \subseteq T \setminus S } \cF_{U} \]
and
\[ U \not\subseteq W \ \ \ \Rightarrow \ \ \ \PP ( \cE_T \wedge \cF_{U}) = 0.\]

Next, let $I = S \cap T$, and define $ \cI$ to be the
event that the process detailed above resorts to revealing potential edges between $v$ and $I$ for more than $ \log^2(n) k $ vertices $ v \not\in S$. 
Note that the probability that a particular vertex $ v \not\in S $ resorts to observing potential edges to $I$ is 
\[ 1-  \frac{ 1-( (1-p)^{\tilde{k}-i} + p(\tilde{k}-i) ( 1- p)^{ \tilde{k}-i-1} )}{ 1 - (1-p)^{\tilde{k}} - p \tilde{k} ( 1- p)^{ \tilde{k}-1}} 
= O \left( \frac{ \log(n) k^2}{ n^2}   \right).\]
Thus, the expected number of vertices that observe edges to $ I$ is $  O(\log(n)^2 k^2/n) = o(k)$. As these events are independent, the Chernoff bound (see Corollary~21.9 in \cite{alanmichal}) then implies
\begin{equation}
\label{eq:I}
\PP( \cI) \le  \exp \{ - \log^2(n) k  \}.
\end{equation}

Now we are ready to put everything together.  Applying the
law of total probability, we have 
\begin{equation*}
\begin{split}
\PP( \cE_T \mid \cE )  & = \sum_{ U \subseteq T \setminus S}
 \left(\PP( \cE_T \mid \cF_{U} \wedge \cI) \cdot \frac{ \PP( \cF_{U} \wedge \cI)}{ \PP( \cE)} +
 \PP( \cE_T \mid \cF_{U} \wedge \overline{\cI}) \cdot \frac{ \PP( \cF_{U} \wedge \overline{ \cI})}{ \PP( \cE)} \right)
\\ 
& \le   \frac{ \PP( \cI)}{ \PP( \cE ) } +   \sum_{ U \subseteq  T \setminus S }
 \PP( \cE_T \mid \cF_{U} \wedge \overline{\cI}) \cdot \frac{ \PP( \cF_{U} )}{ \PP( \cE)} 
\\ 
& =   \frac{ \PP( \cI)}{ \PP( \cE ) } +   \sum_{ U \subseteq  W }
 \PP( \cE_T \mid \cF_{U} \wedge \overline{\cI}) \cdot \frac{ \PP( \cF_{U} )}{ \PP( \cE)} .
\end{split}
\end{equation*}
Now observe that we have
\[
\frac{ \PP( \cF_{U} )}{ \PP( \cE )} 
\le  \left( \frac{ (1-p)^{\tilde{k}-i} + p( \tilde{k}-i) (1-p)^{\tilde{k}-i-1}}{ 1 - (1-p)^{\tilde{k}} - p \tilde{k}( 1-p)^{\tilde{k}-1} }    
\right)^{|U|}  =  \left[ (1+o(1))\frac{ pk^3}{ e^2 n^2} \right] ^{|U|},
\]
and if $ U \subseteq W $ then we have
\begin{equation*}
\begin{split}
\PP( \cE_T \mid \cF_{U} \wedge \overline{\cI})
& \le
\left( \frac{1}{ i } \right)^{|U|}
(1-p)^{ \binom{ \tilde{k}}{2} - \binom{i}{2} - i |U| - (r - |U| - \ell)}  p^{r - |U| - \ell} \\
& \hskip1cm \cdot \left( 1  - (1-p)^{ \tilde{k}} - \tilde{k}p( 1-p)^{\tilde{k}-1}  \right)^{n-2 \tilde{k} +i - \log(n)^2 k} \\
& \le  \frac{ \PP( \cE )}{ p^\ell (1-p)^{i^2/2}} \left( \frac{1}{ ip }  \right)^{|U|} 
\exp\left\{ p \cdot 2ir \right\} \exp\left\{ O\left( \frac{k^3p}{n^2} \left(\tilde{k} + \log(n)^2 k \right) \right) \right\} \\
& \le (1+o(1)) \frac{ \PP( \cE )}{ p^\ell (1-p)^{i^2/2}} \left( \frac{1}{ ip}  \right)^{|U|} 
\exp\left\{ p \cdot 2ir \right\},
\end{split}
\end{equation*}
where we use $ |U| \le \min \{i,r\} $ and $ k^4 p < n^{2 - 2\ep} $.
Thus
\begin{multline*}
\sum_{ U \subseteq  W }
 \PP( \cE_T \mid \cF_{U} \wedge \overline{\cI}) \cdot \frac{ \PP( \cF_{U} )}{ \PP( \cE)}   \\
\le \frac{ (1+o(1)) e^{O(k^3p^2i/n)} \PP( \cE )}{ p^\ell (1-p)^{i^2/2}} \sum_{u=0}^{w} \binom{w}{u} 
 \left[ (1+o(1))\frac{ k^3}{ e^2 i n^2} \right] ^{u} =  \frac{ (1+o(1)) e^{O(k^3p^2i/n)} \PP( \cE )}{ p^\ell (1-p)^{i^2/2}},
\end{multline*}
where we use $ w \le i$ to bound the sum.
We finally note that we have
\[  \frac{ \PP( \cI)}{ \PP( \cE)^2 } \le  \exp \{ - \log^2(n) k + O(\log(n) k) \} = \exp \{ - \Omega\left( \log^2(n)k\right) \} 
= o(1). \]
And the proof of the Lemma is complete.
\end{proof}

\begin{lem} \label{lem:count}
For $\ell =1, \dots, r$ we have
\[ | \bS_\ell| \le | \bS| \left(  \frac{ 3 r^2}{ n^2} \right)^\ell. \]
Furthermore, we have
\[ | \bS_{\ell,i}| \le |\bS| 
\left( \frac{ 3 r^2i^2 }{ \Tilde{k}^4 }\right)^\ell \left( \frac{ e\Tilde{k}^2}{ ni}  \right)^i. \]
\end{lem}

\begin{proof}

We have
\[ | \bS_\ell| \le  \binom{r}{\ell} \binom{n-2\ell}{\Tilde{k}-2\ell} 
\frac{ (\Tilde{k}-2\ell)!}{ (\Tilde{k}-2r)! 2^{r-\ell} (r-\ell)!}. \]
Thus
\begin{equation*}
\frac{ | \bS_\ell| }{ | \bS| } \le \binom{r}{\ell} \frac{ \Tilde{k}(\Tilde{k}-1) \cdots (\Tilde{k} -2\ell+1)}{ n(n-1) \cdots (n-2\ell+1)  }  \cdot
\frac{ 2^\ell \cdot r(r-1) \cdots (r-\ell+1) }{\Tilde{k}( \Tilde{k}-1) \cdot (\Tilde{k} - 2 \ell+1)} 
\le \left( \frac{3 r^2}{ n^2}  \right)^\ell.
\end{equation*}

The calculation for the second part of the Lemma is similar.
\[ | \bS_{\ell,i}| \le  \binom{r}{\ell} \binom{\Tilde{k}-2\ell}{ i -2 \ell} \binom{n-\Tilde{k}}{\Tilde{k}-i} 
\frac{ (\Tilde{k}-2\ell)!}{ (\Tilde{k}-2r)! 2^{r-\ell} (r-\ell)!} \]
Thus
\begin{equation*}
\begin{split}
\frac{ | \bS_{\ell,i}| }{ | \bS| } &\le  \frac{  \binom{r}{\ell}  \binom{\Tilde{k}-2\ell}{ i -2 \ell}\binom{n-\Tilde{k}}{\Tilde{k}-i} \binom{\Tilde{k}}{i}  }{ \binom{n}{i} \binom{n-i}{ \Tilde{k}-i}}  
\cdot \frac{ 2^\ell \cdot r(r-1) \cdots (r-\ell+1) }{\Tilde{k}( \Tilde{k}-1) \cdot (\Tilde{k} - 2 \ell+1)} \\
& \le \binom{r}{\ell} \left( \frac{ e\Tilde{k}^2}{ ni}  \right)^i \left( \frac{i}{ \Tilde{k}} \right)^{2\ell}
\cdot \frac{ 2^\ell \cdot r(r-1) \cdots (r-\ell+1) }{\Tilde{k}( \Tilde{k}-1) \cdot (\Tilde{k} - 2 \ell+1)}  \\
& \le \left( \frac{ 3 r^2i^2 }{ \Tilde{k}^4 }\right)^\ell \left( \frac{ e\Tilde{k}^2}{ ni}  \right)^i .
\end{split}
\end{equation*}

\end{proof}

We are now ready to show that $\sum_{i=0}^{1/(p \log(n))} h_i = o(1)$ using Lemmas \ref{lem:conditioning} and \ref{lem:count}. We will split this into two sub-ranges.

\noindent
{\bf Subcase 1.1}: $i \leq B$

\noindent
In this case, by Lemma \ref{lem:conditioning} we have $\PP(\cE_T \mid \cE) \leq (1+o(1))\frac{\PP(\cE)}{p^{\ell}}$, hence
\begin{equation*}
\begin{split}
\sum_{i=0}^{B}h_i 
& =  \frac{1}{|\bS|} \sum_{\ell = 0}^r \sum_{ (T, m_T) \in \bS_\ell' } \left(  \frac{ \PP( \cE_T \mid \cE)}{ \PP ( \cE) } - 1 \right) \\
& =   \frac{ o( | \bS_0'| )}{ |\bS|} +  \sum_{\ell=1}^r \left( \frac{ 4 r^2}{ p n^2} \right)^\ell   \\
& \le o(1) +  \sum_{\ell=1}^r \left( \frac{1}{ n^{2/3}} \right)^\ell \\
& = o(1).
\end{split}
\end{equation*}

\noindent
{\bf Subcase 1.2}: $B < i < B^2$
\begin{align*}
    h_i 
    & =  \frac{1}{|\bS|} \sum_{\ell=0}^r \sum_{ (T, m_T) \in \bS_{\ell,i}} \left(  \frac{ \PP( \cE_T \mid \cE)}{ \PP ( \cE) } - 1 \right) \\
    &\leq \ds\sum_{\ell=0}^{r}\left(\frac{|\bS_{\ell,i}|}{|\bS |} \right)  \frac{  e^{O(k^3p^2i/n) +o(1)}}{ p^\ell (1-p)^{\binom{i}{2}}} \\&\leq
    \sum_{\ell=0}^{r}\left(\frac{3r^2 i^2}{\Tilde{k}^4}\right)^{\ell}\left(\frac{e \Tilde{k}^2}{ni}\right)^i \left(\frac{2e^{O(k^3p^2i/n)}}{p^{\ell}(1-p)^{\binom{i}{2}}}\right) \\&\leq
    \left(\frac{e^2\Tilde{k}^2}{ni(1-p)^{i/2}}\right)^i \sum_{\ell=0}^{r} \left(\frac{3r^2}{\Tilde{k}^2 p}\right)^{\ell} \\&\leq
    \left(\frac{e^3\Tilde{k}^2 p^{1/2} \ln(n)^{1/2}}{n}\right)^i \sum_{\ell=0}^{r} \left(\frac{3r^2}{\Tilde{k}^2 p}\right)^{\ell}.
\end{align*}
Since $ \frac{ r^2}{ \tilde{k^2}{p}} = O \left( \frac{ k^4p}{n^2}   \right) = o(1)$ and $\frac{\Tilde{k}^2p^{1/2}\log(n)^{1/2}}{n} = o(1)$, we have
\begin{align*}
    \sum_{i=B}^{B^2} h_i = o(1).
\end{align*}

\subsubsection{{\bf Case 2: }$ \frac{1}{ \log(n)p} \leq  i < (1-\ep)\Tilde{k}$}

In this case we can afford to be lenient with our bounds; in particular, we do not need to make use of the condition that
vertices outside $S$ have at least 2 neighbors in $S$. We begin by bounding $ \frac{1}{ \EE( Z_{S,m_S})} \sum_{m_S, m_T} \EE( Z_{S, m_S} Z_{T, m_T}) $ for a fixed pair of sets $S,T$ such that $ |S \cap T| =i$. If there are $\ell$ edges in the intersection, a bound for the number of ways to choose $m_S$ and $m_T$ is $\Tilde{k}^{4r-2\ell}$ (we have $2r-\ell$ edges total; each has less than $\Tilde{k}^2$ choices.) Thus, we have 
\begin{multline*}
\sum_{m_S, m_T} \EE( Z_{S, m_S} Z_{T, m_T})  \le \sum_{\ell=0}^r  \Tilde{k}^{4r-2\ell} p^{2r-\ell} (1-p)^{2\binom{\Tilde{k}}{2}-\binom{i}{2}-2r+\ell}\\
= \Tilde{k}^{4r}p^{2r}(1-p)^{2\binom{\Tilde{k}}{2}-\binom{i}{2}-2r}  \sum_{\ell=0}^r \left( \frac{1-p}{ \Tilde{k}^2 p} \right)^{\ell}
= (1+o(1))\Tilde{k}^{4r}p^{2r}(1-p)^{2\binom{\Tilde{k}}{2}-\binom{i}{2}-2r}   . 
\end{multline*}
Recalling (\ref{eq:maxxy}) we have
\begin{align*}
     \sum_{m_S, m_T}
     \frac{\EE(Z_{S, m_S}Z_{T, m_T})}{\EE(Z_{S,m_S})^2} &\leq (1+o(1))\Tilde{k}^{4r} (1-p)^{-\binom{i}{2}}e^{(2+o(1))\frac{pk^3}{e^2n}} \\&\leq
     (1-p)^{-\binom{i}{2}} e^{ O( \log(n)r+1)}.
\end{align*}
It follows that
\begin{align*}
     h_i & \le \frac{\binom{\Tilde{k}}{i}\binom{n-\Tilde{k}}{\Tilde{k}-i}}{\binom{n}{\Tilde{k}}} (1-p)^{-\binom{i}{2}} e^{ O( \log(n)r+1)} \\
     & \le
     \left(  \frac{ e\Tilde{k}}{i} \right)^i \left( \frac{\Tilde{k}}{n} \right)^i (1-p)^{-\binom{i}{2}} e^{ O( \log(n)r+1)}  \\
     & \le \left[  \frac{e \Tilde{k}^2}{i n}  \left( (1-p)^{-\Tilde{k}/2} \right)^{i/\Tilde{k}}  \right]^i e^{ O( \log(n)r+1)}\\        
     & \le \left[ \frac{e^2 \log(n) p \Tilde{k}^2}{n}  \left( \frac{ne}{\Tilde{k}} \right)^{1-\epsilon}  \right]^i  e^{ O( \log(n)r+1)}  \\ 
     & \le \left[ \frac{ e^{3-\epsilon} \log(n) p \tilde{k}^{1+ \epsilon}}{ n^\epsilon} \right]^i e^{ O( \log(n)r+1)}.
\end{align*}
As $ r  =o(k^{1/2}) = o(i) $, we have
\begin{align*}
    \sum_{i=1/(p \log(n))}^{(1-\ep)\Tilde{k}} h_i = o(1).
\end{align*}

\subsubsection{{\bf Case 3: }$i \geq (1-\ep)\Tilde{k}$}

We calculate $h_i$ as defined in (\ref{h''}). For simplicity, let $j = \tilde{k}-i$. For some fixed $i$ (hence $j$), there are $\binom{n}{\Tilde{k}} \binom{n-\Tilde{k}}{j} \binom{\Tilde{k}}{j}$ ways to choose $S$ and $T$. Now fix sets $S$ and $T$. 

We now consider the number of ways to choose the matchings $m_S$ and $m_T$ that are compatible with the event $ \{ Z_{S,m_S} Z_{T, m_T} =1 \}$. Let $R_I$ be the set of matching edges contained in $S \cap T$, and let $r_I:=|R_I|$ (so $r_I$ is the same as $\ell$ from Cases~1~and~2). Secondly, let $R_O$ be the set of edges in $m_S$ that are not in $R_I$, and let $R'_O$ be the set of edges of $m_T$ that are not in $R_I$, and define $r_O := |R_O| = |R'_O|$ (``$I$" stands for ``inner", and ``$O$" stands for ``outer"). Hence $r = r_O + r_I$. Finally, let $R_S$ be the matching edges in $R_O$ which are contained completely in $S \setminus T$, and let $r_S = |R_S|$. See Figure 1. 

\begin{center} 
    \begin{asy}
size(0,200);
import geometry;
draw(circle((-5,0), 10));
draw(circle((5,0), 10));
draw((-12.5,4)--(-8.5,4));
draw((-12.5,2)--(-8.5,2));
draw((-12.5,0)--(-8.5,0));
draw((-12.5,-2)--(-8.5,-2));
draw((-12.5,-4)--(-8.5,-4));
draw(ellipse((-10.5,0), 3.4, 6),dashed);

draw((12.5,3)--(8.5,3));
draw((12.5,1)--(8.5,1));
draw((12.5,-1)--(8.5,-1));
draw((12.5,-3)--(8.5,-3));
draw(ellipse((10.5,0), 3.2, 5),dashed);

draw((6,0)--(3,0));
draw((5.94,2.92)--(3.06,2.08));
draw((5.94,-2.92)--(3.06,-2.08));
draw(ellipse((4.5,0), 2.3, 5),dashed);

draw((-5.94,2.42)--(-3.06,1.58));
draw((-5.94,-2.42)--(-3.06,-1.58));
draw(ellipse((-4.5,0), 2.3, 4),dashed);

draw((-1,-3)--(1,-3));
draw((-1,-2)--(1,-2));
draw((-1,-1)--(1,-1));
draw((-1,0)--(1,0));
draw((-1,1)--(1,1));
draw((-1,2)--(1,2));
draw((-1,3)--(1,3));
draw(ellipse((0,0),2,5.5),dashed);

label("$S$",(-5,10),N,p=fontsize(18pt));
label("$T$",(5,10),N,p=fontsize(18pt));
label("$R_O$",(-7,-12));
draw((-10.5,-6)--(-7.5,-10.8));
draw((-4.5,-4)--(-7.5,-10.8));

label("$R'_O$",(7.7,-12));
draw((10.5,-5)--(7.5,-10.8));
draw((4.5,-5)--(7.5,-10.8));

label("$R_I$",(0,-12));
draw((0,-5.5)--(0,-10.8));

label("$R_S$",(-14,-12));
draw((-12,-5.5)--(-13.5,-10.8));

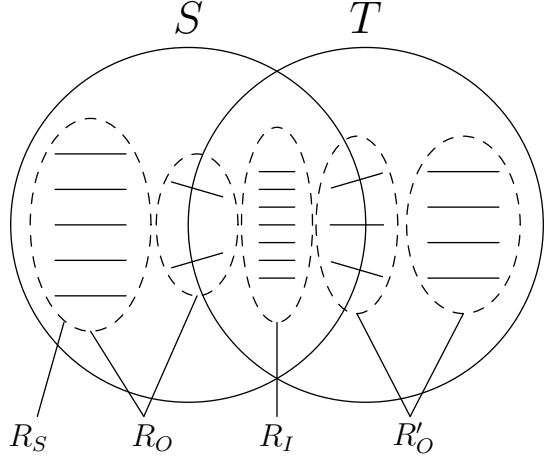
label("{\bf Figure 1.} A sketch of the edges in two intesecting augmented independent sets.", (0,-16));

    \end{asy}
\end{center}

We now count the number of ways to choose edges in $R_O$ if $r_O$ is fixed. To do this, we first note that if $r_S$ is also fixed then the number of ways to choose the edges that comprise $R_S$ is
$$\frac{j!}{(j-2r_S)! r_S! 2^{r_S}}.$$ Then the number of ways to choose the rest of the rest of the edges in $R_O$ is $$\left(\frac{(j-2r_S)!}{(j-r_S-r_O)!}\right)\left(\frac{(\Tilde{k}-j)!}{(\Tilde{k}-j-(r_O-r_S))!}\right) \left(\frac{1}{(r_O-r_S)!}\right),$$ therefore the total number of ways to choose $R_O$ is
\begin{align}
    \sum_{r_S} &\frac{1}{2^{r_S}}\left(\frac{j!}{r_S!(r_O-r_S)!}\right) \left(\frac{(\Tilde{k}-j)!}{(\Tilde{k}-j-(r_O-r_S))!(j-r_S-r_O)! }\right) \nonumber \\&\leq
    \sum_{r_S} \left(\frac{j!}{r_S!(r_O-r_S)!}\right)\left(\frac{(\Tilde{k}-j)^{r_O-r_S}}{(j-r_S-r_O)!}\right) \nonumber \\&\leq
    (\Tilde{k}-j)^{r_O} \sum_{r_S}
    \frac{ j! }{r_S!(r_O-r_S)! (j-r_O)!} \nonumber \\&\leq
    3^j \Tilde{k}^{r_O}.
    \label{rO count}
\end{align}
By a symmetric argument, (\ref{rO count}) is an upper bound for the number of ways to choose the edges in $R'_O$ as well.

Next, note that the number of ways to choose edges in $R_I$ is bounded by $\Tilde{k}^{2r_I}/(2^{r_I}r_I!)$. Putting these together, the number of ways to choose $m_S$ and $m_T$ is at most
\begin{equation}
\label{eq:matchings}
\frac{9^j \Tilde{k}^{2r}}{ 2^{r_I}r_I!} 
\leq
\frac{9^j \Tilde{k}^{2r}  (2r)^{r_O}}{2^{r}r!}
\end{equation}
(assuming that $r^{r_O} = 1$ if $r = r_O = 0$).

Now we estimate the probability of the event $ \{ Z_{S,{m_S}} Z_{T,m_T} = 1 \}$. The
probability the edges within $S$ and within $T$ are chosen accordingly is 
\[ p^{r+r_O}(1-p)^{\binom{\Tilde{k}}{2} + (\Tilde{k}-j)j + \binom{j}{2} -r-r_O}. \] 
Now consider pairs of vertices from $S \backslash T$ to $T \backslash S$. Recall that a vertex $ v \in S \setminus T$ has at least two neighbors in $T$. So, if $ v $ is already incident with an edge that goes into $T$ among the edges already specified, then there is at least one edge from $v$ to $T \backslash S$ (since each from $S \backslash T$ can send at most one edge to $S \cap T$).  Such an edge appears with probability $1 - (1-p)^{j} \le jp$. If a vertex $ v \in S \backslash T$ is not in an edge of $m_S$ that intersects $ S \cap T$ then $v$ is incident with at least two edges going to $T \backslash S$. The probability of these two edges appearing is at most $j^2 p^2$ by the union bound. Since the number of vertices in the first category is at most $r_O$, the total probability of all edges between $S \backslash T$ and $T \backslash S$ being chosen accordingly is at most $(jp)^{2j-r_O}$. Hence, the probability that the edges within $ S \cup T$ appear in accordance with the event $ \{ Z_{S,{m_S}} Z_{T,m_T} = 1 \}$ is at most
\[  p^{r+r_O}(1-p)^{ \binom{\Tilde{k}}{2} + (\Tilde{k}-j)j + \binom{j}{2}-r-r_O}(jp)^{2j-r_O} . \]  
Finally, the probability that all vertices outside $S \cup T$ have at least two neighbors in $S$ is at most $(1-(1-p)^{\Tilde{k}} - \Tilde{k}p(1-p)^{\Tilde{k}-1})^{n-2\Tilde{k}}$. Note that, since
\begin{align*}
    (1-(1-p)^{\Tilde{k}} - \Tilde{k}p(1-p)^{\Tilde{k}-1})^{-\Tilde{k}} &\leq \exp\{ (1+o(1))k^2p(1-p)^k \} \\&=
    \exp\{(1+o(1))k^4p/(e^2n^2)\} \\&=
    1+o(1),
\end{align*}
we have
\[ \PP ( Z_{S,m_S} Z_{T,m_T} = 1) \le (1+o(1)) \PP( Z_{S, m_S} =1 ) p^{r_O}(1-p)^{(\Tilde{k}-j)j + \binom{j}{2} -r_O}(jp)^{2j-r_O}. \]

Multiplying this probability estimate by the estimate for the number of choices for $ m_S$ and $m_T$ given by (\ref{eq:matchings}) and summing over $r_O$ we
see that for a fixed $S$ and $T$ we have
\[  \frac{ 1 }{\PP( Z_{S, m_S} =1 ) }  \sum_{ m_S, m_T} \PP\left(   Z_{S,m_S} Z_{T,m_T} = 1 \right) \]

is at most 
\begin{align*}
    (1+o(1))\sum_{r_O=0}^{\min\{r,j\}} & \frac{9^j \Tilde{k}^{2r} (2r)^{r_O}}{2^r r!} p^{r_O}(1-p)^{ (\Tilde{k}-j)j + \binom{j}{2} -r_O} (jp)^{2j-r_O}
    \\&=
    (1+o(1)) \frac{\Tilde{k}^{2r}}{2^r r!} 
    \left(9p^2j^2 (1-p)^{\Tilde{k}-(j+1)/2}\right)^j \sum_{r_O = 0}^{\min\{r,j\}} \left(\frac{2r}{j(1-p)}\right)^{r_O}
    \\&\leq
    (1+o(1))\frac{\Tilde{k}^{2r}}{2^r r!} 
    \left(9p^2j^2 (1-p)^{\Tilde{k}-(j+1)/2}\right)^j \max\{1,2 (3r)^j\} \\&\leq
    (2+o(1))\frac{\Tilde{k}^{2r}}{2^r r!} 
    \left(30p^2j^2 (1-p)^{\Tilde{k}-j/2 }\max\{1,r\}\right)^j \\&\leq
    (3+o(1)) \frac{\Tilde{k}!}{(\Tilde{k}-2r)!2^r r!} 
    \left(30p^2j^2 (1-p)^{\Tilde{k}-j/2}\max\{1,r\}\right)^j,
\end{align*}
where we use $ r = o( k^{1/2}) $ in the last step (and assume $r^{r_O} = 1$ if $r = r_O = 0$).
Recalling that $ | \bS | = \binom{ n}{\Tilde{k}} \frac{\Tilde{k}!}{(\Tilde{k}-2r)!2^r r!} $ and that the number of choices for the pair of sets $S,T$ is
$ \binom{n}{ \Tilde{k}} \binom{n-\Tilde{k}}{j} \binom{\Tilde{k}}{j}$, we have
\begin{align*}
    h_i \EE(Z) &\leq \frac{1}{ | {\bS}| \PP( Z_{S, m_S} = 1) } \sum_{|S \cap T|=i} \sum_{ m_S, m_T} \PP ( Z_{ S, m_S} Z_{T, m_T} =1 ) \\&\leq
    (3+o(1))\binom{n-\Tilde{k}}{j} \binom{\Tilde{k}}{j} \left(30p^2j^2 (1-p)^{\Tilde{k}-j/2}\max\{1,r\}\right)^j \\&\leq
    (3+o(1))\left(\frac{n\Tilde{k}e^2}{j^2}\right)^j\left(30p^2j^2 (1-p)^{\Tilde{k}-j/2}\max\{1,r\}\right)^j \\&=
    (3+o(1))\left(30e^2n\Tilde{k}p^2(1-p)^{\Tilde{k}-j/2}\max\{1,r\}\right)^j \\&\leq
    (3+o(1)) \left(31e^{\ep}\Tilde{k}^{3-\ep}p^2n^{-1+\ep}\max\{1,r\}\right)^j.
\end{align*} 
As $ r = o( k^{1/2})$ and $ k = O ( \log(n) n^{2/3 - \epsilon})$, we have $\sum_{i=(1-\ep)\Tilde{k}}^{\Tilde{k}-1} h_i = o(1)$, as desired.

\section{Weak anti-concentration of $\a(G_{n,p})$ for $ n^{-1} < p < n^{-2/3}$ } \label{sec:ss}

In this Section we present the argument of Sah and Sawhney that shows that
two-point concentration of the independence number does not extend 
to $p = n^{\gamma}$ if $\gamma \leq -2/3$; that is, we prove Theorem~\ref{SS}. 
We note in passing that this proof is similar to an argument that appears 
in \cite{ak} (see the last item in the `concluding remarks and open problems' 
section of \cite{ak}).
Throughout this section we assume
$ \omega( 1/n ) < p< ( \log(n)/n)^{2/3} $. Recall Theorem~\ref{SS}:
\setcounter{prop}{1}
\begin{thm}[Sah and Sawhney \cite{sahs}] \label{thm:weak}
Let $p = p(n)$ satisfy $ \omega( 1/n ) < p< ( \log(n)/(12n))^{2/3} $  and set $$ \ell =  n^{-1}p^{-3/2} \log(np) / 2.$$ Then there exists 
$ q = q(n)$ such that $p \leq q \leq 2p$ such that $ \alpha( G_{n,q})$ is not
concentrated on $ \ell$ values.
\end{thm}
\setcounter{prop}{10}
\noindent
For $p$ in the specified range we define 
\[ p' = p + n^{-1} \sqrt{p}.\]
We first observe that this choice of $p'$ is close enough to $p$ to ensure
that there is no `separation' of the intervals over which $ G_{n,p}$ and $ G_{n,p'}$
are respectively concentrated. To be precise, we establish the following:

\begin{lem} \label{lem:perturb}
    If $ \omega(n^{-1}) < p < o(1)$ then there is a sequence $ k = k(n)$ such that
    \begin{align*}
        \PP[\a(G_{n,p}) \leq k] > \frac{1}{20} \ \ \ \text{ and } \ \ \  \PP[ \a(G_{n,p'}) \geq k] > \frac{1}{20}
    \end{align*}
    for $n$ sufficiently large.

\end{lem}
\begin{proof}
    First, 
    observe that since the distributions of $e(G_{n,p})$ and $e\left(G_{n,p'}\right)$ are approximately Gaussian with equal variances and with means that are $1/\sqrt{2}$ standard deviations apart, there
    is some value $ m$ such that
    \[  \PP[ e( G_{n,p} ) > m ] > \frac{1}{10} \ \ \ \text{ and } \ \ \ \PP[ e( G_{n,p'} ) < m ] > \frac{1}{10}.  \]
    Now let $k$ be the median value of $ \alpha( G_{n,m} )$. As the addition of edges does not increase the independence number, the Lemma follows.

\end{proof}
Now we prove Theorem~\ref{thm:weak}.
Assume for the sake of contradiction that if $n$ is sufficiently large and $ p \le q \le 2p$
then there is an interval $I$ such that $ |I| \le \ell$ and 
\[ \PP[ \a( G_{n,q}) \in I ] > \frac{49}{50}. \]
Consider the sequence $p_0, p_1, p_2, \dots$, where $p_0 = p$ and $p_{i+1} = p_i'$ for $i \geq 0$, and let $z$ be the  lowest integer such that $p_z \geq 2p$; it is easy to see that $z \leq n\sqrt{p}$.
Let $ I_0 = [a_0,b_0], I_1 = [a_1, b_1], \dots, I_{z-1} = [a_{z-1},b_{z-1}]$ be intervals of $\ell$ values such
that $ \a( G_{n,p_i} )$ lies in $ I_i$ with probability at least $ 49/50$ for each $i$.
By Lemma~\ref{lem:perturb} we have $a_{i} \leq b_{i+1}$ for all $i < z$. This implies 
$ b_{i+1} \ge b_i - \ell$. Iterating this observation gives
\[ b_z \ge b_0 - z \ell \geq b_0 -  \frac{ \log{ np}}{ 2p}. \]

On the other hand, by (\ref{eq:history}), we have 
\[ b_0 - b_z = (1+o(1)) \frac{2 \log(np) }{p} - (1+o(1)) \frac{ 2 \log(2np)}{ 2p} = (1+o(1)) \frac{ \log(np)}{ p}.  \]

This is a contradiction.

\section{Anti-concentration of $\a(G_{n,p})$ for $p = c/n$} \label{sec:anti}

In this Section we establish anti-concentration of the independence number of $ G_{n,p}$ for 
$ p = c/n$ where $c$ is a constant.

We note in passing that for very small $p(n)$ it is easy to show that $\a(G)$ is not narrowly concentrated. Indeed, if $p = n^{\gamma}$ for some $\gamma \in (-2,-3/2)$, then whp no two edges intersect, and hence the independence number is determined by the number of edges that appear. Since the distribution of the number of edges is roughly Poisson with mean $\mu \approx n^{2+\gamma}/2$, then we will not have concentration of $\a(G)$ on fewer than $\Theta(n^{1+\gamma/2})$ values.

We now state our anti-concentration result.

\begin{thm} \label{anti-concentration}
Let $p = c/n$ where $ c>0$ is a constant and let $ k = k(n)$ be an arbitrary sequence of integers.  There exists a constant $L>0$ such that $ \PP( \a(G_{n,p}) = k) \le 
Ln^{-1/2}$ for $n$ sufficiently large. 

\end{thm}
\noindent
Theorem~\ref{anti-concentration} implies that $ \alpha( G_{n, c/n})$ is not concentrated on fewer than $ \sqrt{n}/(2L)$ values. 
Note that this is the best possible anti-concentration result in this regime (up to the constant) as standard martingale methods (such as Hoeffding-Azuma or Talagrand's inequality) show that, for any function $k(n) = \omega(\sqrt{n})$, the independence number is concentrated on some interval of length $k(n)$ (this assertion holds for any value of $p$).

For a fixed graph $G$, let $T (G)$ be the subgraph of $G$ given by all connected components that are trees on $4$ vertices.
Let $t(G)$ be the number of components in $T(G)$. Finally let $\bar{T}(G)$ be the remainder of the graph $G$  (so $ G = T(G) + \overline{T}(G)$). 
Let $p =c/n$ and let $ G = G_{n,p}$. We
set $ \mu = \mu(c)=\frac{2 c^3 e^{-4c}}{ 3}$. Note that
we have
\[ E[t(G)] = \binom{n}{4} \cdot 16 \cdot p^3 (1-p)^{3+4(n-4)}  =(1+o(1)) \frac{16 c^3 e^{-4c}n}{ 24} = (1+o(1)) \mu n . \]
We note that $t(G)$ is a well concentrated random variable.
\begin{lem}
\[ \PP( t(G) < \mu n/2) = O \left( \frac{1}{n}  \right). \]
\end{lem}

\begin{proof}
Here we use standard second moment method once again. For vertex set $S$ with $|S| = 4$, let $W_S$ be the random variable that indicates if the graph induced on $S$ is an isolated tree. Set $ W =t(G)$ and note that $ t(G)$ is the sum of all $\binom{n}{4}$ indicators.  
By Chebychev, we have
\begin{align*}
     \PP\left(W < \frac{\mu n}{2}\right) \leq \PP\left(|W - \EE(W)| \geq \frac{\EE(W)}{3}\right) \leq \frac{9\text{Var}(W)}{\EE(W)^2}.
\end{align*}
Hence
\begin{align*}
     \PP\left(W < \frac{ \mu n}{2}\right) \leq \frac{9\text{Var}(W)}{\EE(W)^2} < 9\left(\frac{1}{\EE(W)} + \max_{S,T}\left\{\frac{\EE(W_S W_{T}) - \EE(W_S)\EE(W_{T})}{\EE(W_S)^2}\right\}\right).
 \end{align*}
Note that the only case where we get a positive covariance is where $S \cap T = \emptyset$, and in this case, we have
\begin{align*}
     \frac{\EE(W_S W_{T}) - \EE(W_S)\EE(W_{T})}{\EE(W_S)^2} = (1-p)^{-16}-1 = \frac{16c+o(1)}{n}.
 \end{align*}
Therefore
\begin{align*}
\PP\left( t(G) < \frac{ \mu n}{2}\right) < \frac{9 +o(1)}{ \mu n} + \frac{144 c+o(1)}{n} =O\left(\frac{1}{n} \right).
\end{align*}
\end{proof}

With this observation in hand, we are ready to proceed to the proof of Theorem~\ref{anti-concentration}. Let $k$ be a fixed integer. 
Let $ \cT$ be the collection of all graphs $H$ on $v_H \le n$ vertices such that $ 4 \mid n - v_H$ and $H$ has no connected component that is a tree on 4 vertices.
Let $ \cT'$ be the collection of graphs in $H \in \cT$ such that $ v_H \le n - 2 \mu n$ vertices.
Applying the law of total
probability we have
\begin{multline*}
\PP( \alpha( G) = k) = \sum_{H \in \cT} \PP( \alpha( G) = k \mid \overline{T}( G) = H) \PP( \overline{T}(G) = H) \\
\le \PP (t(G) \le \mu n/2) +  \sum_{H \in \cT'} \PP( \alpha( G) = k \mid \overline{T}( G) = H) \PP( \overline{T}(G) = H).
\end{multline*}
So, it suffices to bound the conditional probability
\[  \PP( \alpha( G) = k \mid \overline{T}( G) = H) \]
where $H$ is a graph in $ \cT'$. Under this conditioning the graph $G$ is $H$ together with a forest of $ (n - v_H)/4$ trees on $4$ vertices. Furthermore, the number of trees
in this forest which are stars on 4 vertices is a binomial random variable with $ (n - v_H)/4$  trials with probability of success $ 1/4$ on each trial. Let $B$ be this random variable,
and note that if $G$ is drawn from the conditional probability space on the collection of graphs with $\overline{T}( G) = H $ then we have
\[ \alpha(G) = \alpha(H) + (n - v_H)/2 + B. \]
Therefore, we have
\[ \PP( \alpha( G) = k \mid \overline{T}( G) = H) = \PP( B = k -\alpha(H) - (n - v_H)/2 \mid  \overline{T}( G) = H).  \]
But $B$ is a binomial random variable with a linear number of trials, and therefore the probability that $B$ is
any particular value is at most $ L/\sqrt{n}$ for some constant $L$.

\section{Conclusion} \label{sec:concl}

In this Section we present some remarks and open questions 
related to the extent of concentration of $ \alpha( G_{n,p})$ for $ p < n^{-2/3}$. We write
$ f = \tilde{O}(g)$ if $ f/g$ is bounded above by a polylogarithmic function
and 
$ f = \tilde{\Omega}(g)$ if $ f/g$ is bounded below by a polylogarithmic function.
\begin{itemize}

\item The precise extent of concentration of $ \a( G_{n,p})$ for $ p = n^{\gamma}$ with
$ -1 < \gamma < -2/3$ remains an interesting open question. Note that Theorem~\ref{SS} gives a lower bound of $ \tilde{\Omega}( n^{-1- 3\gamma/2})$ on the extent of concentration in this range. Furthermore, Theorems~\ref{MainLoose}~and~\ref{anti-concentration} show that this lower bound 
is of the correct order at either end of the interval. So it is natural to ask if 
this lower bound is actually tight throughout this range.
\begin{quest}
Suppose $ p = n^{\gamma}$ where $\gamma$ is a constant such that $ -1 < \gamma < -2/3$. Is $ \a( G_{n,p} )$ concentrated on $ \tilde{O}\left( n^{-1 - 3 \gamma/2}   \right)$ values?
\end{quest}

\item We do not see a way to adapt the argument given in Section~\ref{sec:ss} to prove
anti-concentration of $ \a(G_{n,m})$ for $m \le n^{4/3}$. We conjecture that
we have such anti-concentration.
\begin{conj}\label{conj:edges}
Suppose $ m = n^{ \eta}$ where $\eta$ is a constant such that $ 1 < \eta < 4/3$. Then $ \a(G_{n,m})$ is not concentrated
on $ n^{ 2 - 3\eta/2 - \epsilon} $ values.
\end{conj}
At the moment we have no lower bound on the extent of concentration $ \a( G_{n,m})$ in this regime.

It seems natural to suspect that Conjecture~\ref{conj:edges} is closely related to the questions
of proving a stronger, uniform 
form of anti-concentration of $ \a( G_{n,p})$ for $ p = n^{\gamma}$ with $ -1 < \gamma < -2/3$. 
We believe that the following strengthening of Theorem~\ref{SS} holds.
\begin{conj}\label{conj:strong}
If $ p = n^{\gamma}$ where $ \gamma$ is a constant such that $-1 < \gamma < -2/3$ then we have
\[ \PP\left( \a( G_{n,p}) = k \right) = \tilde{O} \left( n^{1 + 3 \gamma/2}  \right) \]
for any sequence $k = k(n)$.
\end{conj}
\noindent
Note that we proved Conjecture~\ref{conj:strong} for $ \gamma =-1$ (i.e. $ p = 1/n$) in Section~\ref{sec:anti}.

\item We conclude by noting that the first moment argument using 
augmented independent sets given in Section~\ref{sec:firstmoment} can be
adapted in straightforward way to give the following 
upper bound on the independence number of $ G_{n,p}$ for $ p \le n^{-2/3}$. 
\begin{thm}
If $ \omega( 1/n) < p \le n^{-2/3}$ then we have
\[ \alpha (G_{n,p}) \le k_x - \Omega\left( \frac{ \log(np)^2 }{ p^2n} \right) \]
with high probability, where $k_x$ is the largest integer such that 
$ \binom{n}{k_x}(1-p)^{\binom{k_x}{2}} > 1$. 
\end{thm}
This improvement can be extended to give a small linear improvement on the best known upper
bound on $ \alpha( G_{n,p} )$ where $ p = c/n$ and $c$ is a constant. In fact, a first moment
argument using maximal independent sets is sufficient to give a linear improvement on the
upper bound in this context. As
the precise magnitude of the optimal improvement requires a lengthy calculation, we
exclude this in the interest of brevity.

\end{itemize}

\noindent
{\bf Acknowledgement.} We thank the anonymous referee for useful comments and suggestions.

\pagebreak
\appendix

\noindent{\large \bf Appendix: Two-point concentration of $ \a(G_{n,p})$ for $n^{-1/2+ \epsilon} < p< 1/ \log(n)^2$} 

\noindent
Here we show that the second moment method applied to the random variable $Y_k$, which is the number of
maximal independent sets of size $k$, suffices to establish two point concentration of $ \a(G_{n,p}) $ for $n^{-1/2+ \epsilon} < p< 1/ \log(n)^2$.

Recall that $X_k$ is the random variable which counts the number of independent sets of size $k$ in $ G_{n,p}$, and we defined $k_x$ to be the largest integer such that 
\begin{align}
    \EE(X_k) > n^{2 \epsilon} \label{def of kx}.
\end{align}
Set $ X = X_{k_x}$ and $ Y = Y_{x_k}$.

\begin{thm} \label{maximal}
Consider $p = p(n)$ such that
$ n^{-1/2 + \epsilon} < p< 1/ \log(n)^2$. Then $\a(G_{n,p}) \in \{k_x,k_x+1\}$ whp. 
\end{thm}

\begin{proof}

Set $k = k_x$. First, we show by the standard first moment method that no independent set of size $k+2$ appear whp:
\begin{align*}
\PP( \alpha( G_{n,p}) \ge k +2)
    &\leq \binom{n}{k +2}(1-p)^{\binom{k+2}{2}} \\&= \binom{n}{k+1}(1-p)^{\binom{k+1}{2}}\Bigg(\frac{n-k-1}{k+2}(1-p)^{k+1} \Bigg) \\&=
    (1+o(1))\binom{n}{k+1}(1-p)^{\binom{k+1}{2}} \frac{n}{k} \left( \frac{k}{ne} \right)^2 \\& =
    o \left( \frac{k}{n^{1/2}} \right) \\&=
    o(1),
\end{align*}
where we use the fact that $ k = o(n^{1/2}) $ when $p > n^{-1/2+\ep}$ in the last line.

Next, we show that an independent set of size $k$ exists whp. Here we use the second moment method, but the variance of $X$ itself is too large. So we work with the random variable $Y$ which counts the number of {\it maximal} independent sets of size $k$. The first moment calculation is straightforward, using (\ref{first max factor}):
\begin{align}
    \EE(Y) &= \binom{n}{k}(1-p)^{\binom{k}{2}} (1-(1-p)^{k})^{n-k} \nonumber \\&=
    \binom{n}{k}(1-p)^{\binom{k}{2}}e^{- O( k^2/n)} \label{extra factor}\\&=
    \EE(X)(1-o(1)) \label{first moment modified} \\&= \omega(1). \nonumber
\end{align}

By Chebyshev's Inequality, all that is left to show is that $\frac{Var(Y)}{\EE(Y)^2} = o(1)$. We write $Y$ as the sum of indicator variables $Y_S$, over all sets $S$ such that $|K| = k$, of the random variable $Y_S$ which is the indicator for the event that $S$ is a maximal independent vertex set. 
We define indicator variables $X_S$ for not-necessarily-maximal independent sets analogously. 
Then we have
\begin{align*}
    \text{Var}(Y) = \sum_S \text{Var}(Y_S)+ \sum_{S \neq T} \text{Cov}(Y_S,Y_{T}) \leq \EE(Y)+ \sum_{S \neq T} \EE(Y_S Y_{T}) - \EE(Y_S) \EE(Y_{T})
\end{align*}
Next, define, for all relevant $i$:
\begin{align*}
    &f_i = \frac{1}{\EE(X)^2}\sum_{|S \cap T| = i} \EE(X_S X_{T}) = \frac{\binom{k}{i}\binom{n-k}{k-i}}{\binom{n}{k}} (1-p)^{-\binom{i}{2}}, \\&
    g_i = \frac{1}{\EE(Y)^2}\sum_{|S \cap T| = i} \EE(Y_S Y_{T}) - \EE(Y_S) \EE(Y_{T}), \text{ and} \\&
    \kappa_i = \frac{f_{i+1}}{f_i}.
\end{align*}
Since $\EE(Y) = \omega(1)$, our goal is to show that 
\begin{align*}
    \sum_{i=0}^{k-1} g_i = o(1).
\end{align*}

We consider four ranges of $i$: $i = 0, 1 \leq i < \frac{\ep}{2} k_x$, $\frac{\ep}{2} k_x \leq i \leq (1-\ep)k_x$, and $i > (1-\ep) k_x$ (where $\epsilon$ is the constant in the statement of the
Theorem). For each range, we show that the sum of variables $g_i$ within that range is $o(1)$. For the second and third ranges, we will actually work with variables $f_i$ instead, noting that, since $Y_S \leq X_S$ and by similar reasoning to (\ref{first moment modified}), $g_i \leq f_i(1+o(1))$. 

First, consider $i = 0$. Say we fix two disjoint vertex sets $S$ and $T$ of size $k$; we have $$g_0 = \frac{\EE(Y_S Y_T)}{\EE(Y_S)^2} - 1.$$

Define $Y_S'$ to be the indicator variable of the event that $S$ is an independent set and that every vertex not in $S \cup T$ has a neighbor in $S$; likewise, define $Y_T'$ to be the indicator variable of the event that $T$ is an independent set and that every vertex not in $S \cup T$ has a neighbor in $T$. Finally, define $I$ to be the indicator variable of the event that every vertex in $S$ has a neighbor in $T$ and vice versa. Since $Y_S \geq Y_S' I, Y_S Y_T = Y_S' Y_T' I$, and $Y_S', Y_T', $ and $I$ are independent, we have

$$\frac{\EE(Y_S Y_T)}{\EE(Y_S)^2} \leq \frac{1}{\EE(I)}.$$

Since the probability of all vertices in $S$ having a neighbor in $T$ is $(1 - (1-p)^k)^{k} = 1 - o(1)$, then by union bound $\EE(I) = 1 - o(1)$, hence

$$g_0 \leq \frac{1}{\EE(I)} - 1 = o(1).$$

Next, consider $1 \leq i < \frac{\ep}{2} k$. Using the explicit formula for $f_i$ we have
\begin{multline*}
    \kappa_i = \frac{(k-i)^2(1-p)^{-i}}{(i+1)(n-2k+i+1)} <
    \frac{ k^2 (1-p)^{-\ep k/2}}{n} = (1+o(1)) \frac{ k^2}{n} \left( \frac{ ne}{k} \right)^{\epsilon}
    = O \left( \frac{ \log(n)^2}{ n^{\epsilon}}  \right)
    = o(1).
\end{multline*}
As $ f_{1} = O( k^2/n) = o(1)$ then $\sum_{i=1}^{\ep k/2} f_i$ is a geometric sum with leading term and ratio in $o(1)$.

Next, consider the range $\frac{\epsilon}{2} k \leq i \leq (1-\epsilon) k$. We have
\begin{equation}
    \begin{split}
    f_i & \le 2^{k} \left( \frac{k}{n} \right)^i (1-p)^{-\binom{i}{2}} \\
    & \le \left[  2^{k/i}  \frac{k}{n}  \left( (1-p)^{-k/2} \right)^{i/k}  \right]^i \\        
    & \le \left[ (1+o(1))  2^{ 2/ \epsilon }  \frac{k}{n}  \left( \frac{ne}{k} \right)^{1-\epsilon}  \right]^i \\        
    & \le \left[ (1+o(1))  2^{ 2/ \epsilon } e^{1-\epsilon}  \left( \frac{k}{n}  \right)^\epsilon  \right]^i. \\        
    \end{split} \label{middle}
\end{equation}  
Hence the sum of variables $f_i$ over $\frac{\ep}{2} k_x \leq i \leq (1-\ep)k_x$ is $o(1)$.

The range where $i > (1-\ep) k$ is where we use the variables $g_i$ instead of $f_i$. Given two vertex sets $S$ and $T$ with intersection size $i$, and given that they are both independent, consider the probability that both are maximal independent sets. This probability is bounded above by the probability that each vertex in $S \backslash T$ is adjacent to at least one vertex in $T \backslash S$, which is equal to $(1-(1-p)^{k-i})^{k-i} \leq ((k-i)p)^{k-i}$. Therefore we have
\begin{align*}
    g_i &\leq (1+o(1)) \frac{\binom{k}{i}\binom{n-k}{k-i}}{\binom{n}{k}} (1-p)^{-\binom{i}{2}} ((k-i)p)^{k-i}  \\&= \frac{1+o(1)}{\EE[X]}
    \Bigg(\binom{k}{k-i} (k-i)^{k-i}\Bigg) \Bigg(\binom{n-k}{k-i}(1-p)^{\binom{k}{2}-\binom{i}{2}}p^{k-i} \Bigg)  \\&< \frac{1+o(1)}{\EE[X]}
    (ke)^{k-i} n^{k-i} (1-p)^{(1-\ep)k(k-i)}p^{k-i} \\&< \frac{1}{\EE[X]}
    \Bigg((1+o(1)) e k n p \left( \frac{ k}{ ne} \right)^{2(1 -\epsilon)} \Bigg)^{k-i} \\&< \frac{1}{\EE[X]}
    \Bigg( 22  \log(n)^3 n^{-\epsilon + 2\epsilon^2} \Bigg)^{k-i},
\end{align*}
so $\sum_{i>(1-\ep) k} g_i = o(1)$, as desired.
\end{proof}

\end{document}